\documentclass[a4paper,12pt,leqno]{article}
\usepackage{latexsym}
\usepackage[all]{xy}

\usepackage{amssymb} 
\usepackage{amsmath} 

\usepackage{tikz}
\usetikzlibrary{arrows,cd}  
\usepackage{amsthm}

\def\P{{\mathbf{P}}}
\def\Z{{\mathbb{Z}}}
\def\Q{{\mathbb{Q}}}
\def\K{{\mathbb{K}}}
\def\CC{{\mathbb{C}}}
\def\R{{\mathbb{R}}}
\def\C{{\mathbb{C}}}
\def\A{{\mathcal{A}}}
\def\B{{\mathcal{B}}}
\def\C{{\mathcal{C}}}

\DeclareMathOperator{\codim}{codim}

\DeclareMathOperator{\coker}{coker}
\DeclareMathOperator{\Der}{Der}

\DeclareMathOperator{\pd}{pd}

\DeclareMathOperator{\res}{res}

\DeclareMathOperator{\POexp}{POexp}

\numberwithin{equation}{section}

\newcommand{\owari}{\hfill$\square$}

\newtheorem{theorem}{Theorem}[section]
\newtheorem{prop}[theorem]{Proposition}
\newtheorem{cor}[theorem]{Corollary}

\newtheorem{define}[theorem]{Definition}

\newtheorem{problem}[theorem]{Problem}
\newtheorem{conj}[theorem]{Conjecture}
\theoremstyle{remark}
\newtheorem{rem}[theorem]{Remark}
\newtheorem{example}[theorem]{Example}

\title{Addition-deletion theorems for the Solomon-Terao polynomials and 
$B$-sequences of hyperplane arrangements}
\author{Takuro Abe}
\date{\today}

\begin{document}

\maketitle

\begin{abstract}
We prove the addition-deletion theorems for the Solomon-Terao polynomials, which have two 
important specializations. Namely, one is to the characteristic polynomials of hyperplane arangements, and the other to the Poincar\`{e} polynomials of the regular nilpotent Hessenberg varieties. One of the main tools to show them is the 
free surjection theorem which confirms the right exactness of several important 
exact sequences among logarithmic modules. 
Moreover, we introduce a generalized polynomial $B$-theory to the higher order 
logarithmic modules, whose origin was due to Terao.
 \end{abstract}

\section{Introduction}
Let $\K$ be a field, $V=\K^\ell$ and $S:=\mbox{Sym}^*(V^*) =\K[x_1,\ldots,x_\ell]$, where $x_1,\ldots,x_\ell$ form a basis for $V^*$. Let $\Der S:=\oplus_{i=1}^\ell S\partial_{x_i}$ be the module of 
$\K$-linear $S$-derivations, and 
$\Omega^p_V:=\sum_{1\le i_1 < \cdots <i_p \le \ell} 
S dx_{i_1}\cdots dx_{i_p}$ be the 
\textbf{module of regular $p$-forms}. Let $\A$ be an \textbf{arrangement} of hyperplanes in $V$, i.e., a finite set of linear hyperplanes in $V$. For each $H \in \A$ let $\alpha_H \in V^*$ be its defining linear form, and let $Q(\A):=\prod_{H\in \A} \alpha_H$. Now recall the following algebraic objects of hyperplane arrangements:

\begin{define}
For $0 \le p \le \ell$, define 
the \textbf{logarithmic derivation module of order $p$} by 
\begin{eqnarray*}
D^p(\A):&=&\{\theta \in \Der^p S:=\wedge^p \Der S \mid \\
&\ &
\theta(\alpha_H,f_2,\ldots,f_p) \in S \alpha_H\ (\forall H \in \A,\ \forall f_2,\ldots,f_p \in S)\},
\end{eqnarray*}
and the \textbf{logarithmic module of the differential $p$-forms} 
by 
$$
\Omega^p(\A):=\{\omega \in \displaystyle \frac{1}{Q(\A)} \Omega^p_V \mid 
Q(\A) \omega \wedge d\alpha_H \in \alpha_H \Omega^{p+1}_V\ (\forall H \in \A)\}.
$$
\label{logmodules}
\end{define}

It is known that $D^p(\A)^*\simeq \Omega^p(\A)$ and $\Omega^p(\A)^* \simeq D^p(\A)$ 
(see Theorem 4.75 in \cite{OT} for example), thus they are $S$-graded reflexive modules. However, 
these are very difficult to study when $p>1$. So many researches on these logarithmic modules are for 
the cases of $p=1$.
In particular, the freeness of $\A$ is described in terms of $D(\A)=D^1(\A)$ (just called the 
\textbf{logarithmic derivation module}). Namely, 
we say that $\A$ is \textbf{free} with $\exp(\A)=(d_1,\ldots,d_\ell)$ if $$
D(\A) \simeq \oplus_{i=1}^\ell S[-d_i].
$$
The freeness has played one of the central roles in the study of hyperplane arrangements. 
By using logarithmic derivation modules, Solomon and Terao introduced the following polynomial.

\begin{define}[\cite{ST}, \cite{AMMN}]
The \textbf{Solomon-Terao polynomial} 
$\Psi(\A;x,t)$ is defined by 
$$
\Psi(\A;x,t):=
\sum_{p=0}^\ell \mbox{Hilb}(D^p(\A);x)(t(x-1)-1)^p.
$$
\label{ST}
\end{define}

Clearly $\Psi(\A;x,t)$ is a series, but 
the following surprising result was proved in \cite{ST}.

\begin{theorem}[Theorem 1.2 and 
Proposition 5.3, \cite{ST}]
It holds that 
$$\Psi(\A;x,t) \in \Q[x,t].$$ Moreover, 
$$
(-1)^\ell \Psi(\A;1,t)=\chi(\A;t).
$$
Here $\chi(\A;t)$ is the characteristic polynomial of $\A$, see 
Definition \ref{combin}. 
\label{STformula}
\end{theorem}

So $D^p(\A)$ knows comibinatorics and topology of $\A$, see 
Definition \ref{combin} too. Then it is natural to ask the other specialization, i.e., can we get invariants when some value 
is substituted to $t$? The answer is yes, i.e., 
when $t=-1$, 
the \textbf{reduced Solomon-Terao polynomial} $\Psi(\A;x,-1)$ knows topology of the regular nilpotent 
Hessenberg varieties as follows:

\begin{theorem}[Theorem 1.3 and Corollary 7.3, \cite{AHMMS}]
$\Psi(\A;x,-1)$ is essentially the same as the topological Poincar\`{e} 
polynomial of the regular nilpotent Hessenberg variety $\mbox{Hess}(I)$ 
when $\A$ is the ideal arrangement corresponding to the ideal $I$ in some 
positive system.
\label{Hess}
\end{theorem}

So $\Psi(\A;x,t)$ has two nice specializations. Thus it is natural to ask what $\Psi(\A;x,t)$ itself is, and also what is $\Psi(\A;x,-1)$ 
for other arrangements. However, the research of $\Psi(\A;x,t)$ does not exist except for the two papers \cite{ST} and \cite{AMMN}. In fact, we do not know $\deg \Psi(\A;x,-1)$ in general, see Conjecture \ref{degree}. 

The reason is the hardness to compute $\Psi(\A;x,t)$. In fact the only case they are computed is when $\A$ is free, see Theorem \ref{STfree}. That is based on the difficulty of the structure of $D^p(\A)$ when $2 \le p \le \ell-1$. Hence there 
are no formula like the deletion-restriction that holds for the characteristic polynomials 
as follows:
$$
\chi(\A;t)=\chi(\A';t)-\chi(\A^H;t).
$$
Here $\A^H:=\{ L \cap H \mid 
L \in \A \setminus \{H\}\}$
is the \textbf{restriction} of $\A$ onto $H \in \A$. Since the definition of the Solomon-Terao polynomials are algebraic, in general 
such a formula does not hold. However, there 
is a famous exact sequence called the 
\textbf{Euler exact sequence} as follows:
\begin{equation}
0 \rightarrow D^p(\A') 
\stackrel{\cdot \alpha_H}{\rightarrow} 
D^p(\A) \stackrel{\rho^p}{\rightarrow} 
D^p(\A^H).
\label{euler}
\end{equation}
Here 
$\A':=\A \setminus \{H\}$, and 
$\rho(\theta)(\overline{f}):=\overline{\theta(f)}$ for $f \in S$, where $\overline{f}$ 
denotes the canonical image of $f \in S$ in $\overline{S}=S/\alpha_H S$. 
In general (\ref{euler}) is not right exact. However, in this article, under a certain condition, we give
the addition-deletion formula  for the Solomon-Terao polynomials by 
showing the right exactness of the Euler and the other exact sequences proved in this article later. 

To state it, let us introduce one more terminology.  For $X \in L(\A):=\{ \cap_{H \in \B} H \mid \B \subset \A\}$, the 
\textbf{localization of $\A$ at $X$} is defined by 
$$\A_X:=\{H \in \A 
\mid X \subset H\}.$$
Let $\A$ be an arrangement, $H \in \A$ and let $\B$ be an arrangement in $H$. 
When we have two functors $F,G$ from the category of hyperplane arrangements to the
$S$-modules, we say that a map $f:F(\A) \rightarrow G(\B)$ is 
\textbf{surjective in codimension $k$ along $H$} if the induced map 
$$
f_X:F(\A_X) \rightarrow G(\B_X)
$$
is surjective for all $X \in L(\B)$ with $H \supset X$ and $\codim_V X=k$.
Now let us state the addition theorem as follows:

\begin{theorem}[Addition theorem for Solomon-Terao polynomials]
Let $H \in \A$ and let $\A':=\A\setminus \{H\}$. If 
$\rho^H$ in (\ref{euler}) is surjective in codimension $k+2$ along $H$ and 
$\pd_S D(\A')=k< \ell-2$ for all $0 \le p \le \ell$, 
then 
$$
\Psi(\A;x,t)=
x\Psi(\A';x,t)+\Psi(\A^H;x,t)
$$
and 
$$
\Psi(\A;x,-1)=
x\Psi(\A';x,-1)+\Psi(\A^H;x,-1).
$$
\label{STminus1dual}
\end{theorem}

The conditions in Theorem \ref{STminus1dual} may seem complicated, but for example, 
if $\A'$ is free, the these are automatically satisfied (see Corollary \ref{STexplicit}). 
Next let us 
consider the deletion theorem. For that purpose, for an arrangement $\A$. $H \in \A$ and $\A':=\A
\setminus \{H\}$, let us define maps 
$\overline{\partial}^p:D^p(\A') \rightarrow \wedge \Der^{p-1} \overline{S}$ by 
$$\overline{\partial^p}(\theta)(\overline{f}_2,\ldots,\overline{f}_p):=
\overline{\theta(\alpha_H,f_2,\ldots,f_p)}
$$
for $\theta \in D(\A')$.
Then we obtain the following new exact sequence between 
$D(\A)$ and $D(\A')$. 

\begin{theorem}[$B$-sequence]
There is an exact sequence 
\begin{equation}
0 \rightarrow D^p(\A) \rightarrow D^p(\A') \stackrel{\overline{\partial^p}}{\rightarrow} D^{p-1}(\A^H)\overline{B}.
\label{seB}
\end{equation}
Here $0 \le p \le \ell$, $\overline{B}:=
\displaystyle \frac{\overline{Q(\A')}}{Q(\A^H)} \in \overline{S}$. 
We call the exact sequence (\ref{seB}) the \textbf{$B$-sequence}. 
\label{higherB}
\end{theorem}

Theorem \ref{higherB} answers two questions on these modules, i.e., (1) what is a good cokernel and map of the canonical inclusion $D^p(\A) \subset D^p(\A')$, and (2) there are any polynomial $B$-theory for $D^p(\A)$ as when $p=1$ due to Terao. When $p=1$, we have Terao's polynomial $B$-theory, see Theorem \ref{B}. 
Since $B$ in Theorem \ref{higherB} is the same as Terao's polynomial $B$, we can 
regard Theorem \ref{seB} as a higher version of Theorem \ref{B}. Now by using the $B$-sequence, we can show the deletion theorem for the Solomon-Terao polynomials.

\begin{theorem}[Deletion theorem for Solomon-Terao polynomials]
Let $H \in \A,\ \A':=\A \setminus \{H\}$ and $
d:=|\A'|-|\A^H|$. If $\overline{\partial}^p$ in Theorem \ref{higherB} are surjective in codimension $k+2$ along $H$ and 
$\pd_S D(\A)=k< \ell-2$ for all $0 \le p \le \ell$, then 
$$
\Psi(\A';x,t)=
\Psi(\A;x,t)
+x^d(t(x-1)-1)\Psi(\A^H;x,t).
$$
and 
$$
\Psi(\A';x,-1)=
\Psi(\A;x,-1)-x^{d+1}
\Psi(\A^H;x,-1).
$$
\label{STminus1}
\end{theorem}

As Theorem \ref{STminus1dual}, the conditions in Theorem \ref{STminus1} are satisfied when $\A$ is free 
(see Corollary \ref{STexplicit}). In these cases, we can show explicit formula as follows:

\begin{cor}
Let $H \in \A,\ \A':=\A \setminus \{H\},\ 
d:=|\A'|-|\A^H|$.

(1)\,\,
Assume that $\A'$ and $\A^H$ are both free with 
$\exp(\A')=(1,d_2,\ldots,d_\ell)$ and 
$\exp(\A^H)=(1,e_2,\ldots,e_{\ell-1})$. Then 
$$
\Psi(\A;x,t)=x\prod_{i=1}^\ell (-tx^{d_i}+1+x+\cdots+x^{d_i-1})
+\prod_{i=1}^{\ell-1}(-tx^{e_i}+1+x+\cdots+x^{e_i-1})
$$
and 
$$
\Psi(\A;x,-1)=x\prod_{i=1}^\ell (1+x+\cdots+x^{d_i})
+\prod_{i=1}^{\ell-1}(1+x+\cdots+x^{e_i}).
$$

(2)\,\,
If $\A$ and $\A^H$ are free with 
$\exp(\A)=(1,d_2,\ldots,d_\ell)$ 
and 
$\exp(\A^H)=(1,e_2,\ldots,e_{\ell-1})$, then 
\begin{eqnarray*}
\Psi(\A';x,t)&=&
\prod_{i=1}^\ell 
(-tx^{d_i}+\sum_{i=0}^{d_i-1} x^{i})
+x^d(t(x-1)-1)\prod_{i=1}^{\ell -1}
(-tx^{e_i}+\sum_{i=0}^{e_i-1} x^{i})
\end{eqnarray*}
and 
\begin{eqnarray*}
\Psi(\A';x,-1)&=&
\prod_{i=1}^\ell 
(1+x+\cdots+x^{d_i})-x^{d+1}\prod_{i=1}^{\ell -1}
(1+x+\cdots+x^{e_i}).
\end{eqnarray*}

\label{STexplicit}
\end{cor}

\begin{rem}
Theorems \ref{STminus1dual} and \ref{STminus1} give the deletion-restriction formula when $x=1$ multiplied with $(-1)^\ell$ and Theorem \ref{ST}. So 
it is mysterious that the higher structure, $\Psi(\A;x,t)$ of $\chi(\A;t)$ has at least two types of the addition-deletion theorems 
(Theorems \ref{STminus1dual} and \ref{STminus1}). 
So these two theorems say that if both conditions in them are satisfied, the results are the same, i.e., 
$$
\Psi(\A;x,t)=\Psi(\A';x,t)-x^d(t(x-1)-1)\Psi(\A^H;x,t)=
x\Psi(\A';x,t)+\Psi(\A^H;x,t).
$$
 For example, if both $\A$ and $\A'$ are free, then surprisingly Theorems \ref{STminus1dual} and \ref{STminus1} give the same formula. This is because the proofs of them are highly algebraic, not combinatorial, and that is the reason why we call them the addition-deletion theorem, not the deletion-restriction. 
\end{rem}

Theorem \ref{higherB} has a lot of corollaries. First we can give another proof of Theorem \ref{B} as the case $p=1$. The most interesting corollary is the following SPOG-type theorem for higher derivations (see 
Definition \ref{SPOG} and Theorem \ref{SPOGthm} for details).

\begin{theorem}
Assume that $D^p(\A)$ is free. Then for any $H \in \A$, if 
$D^p(\A')$ is not free for $\A':=
\A \setminus \{H\}$, then 
$$
D^p(\A')=D^p(\A)+\sum_{i=1}^s S \varphi_i,
$$
where $\overline{\partial^p}(\varphi_1),\ldots,\overline{\partial^p}(\varphi_s)$ form a minimal set of generators for $D^{p-1}(\A^H)\overline{B}$. Moreover, 
$$
\pd_S D^p(\A') =\pd_{\overline{S}} D^{p-1}(\A^H)+1.
$$
In particular, if $D^p(\A)$ and $D^p(\A^H)$ are both free, then $\pd_S D^p(\A \setminus \{H\}) \le 1$. 
\label{SPOGhigher}
\end{theorem}

The case when $p=1$ in Theorem \ref{SPOGhigher} corresponds to the 
original SPOG theorem (see Definition \ref{SPOG} and Theorem \ref{SPOGthm}). 
Also by these results we can determine the structure of $D^p(\A')$ when $\A$ is free as follows:

\begin{cor}
Assume that $D^p(\A)$ is free with $\exp^p(\A)=\{d_I\}_{I \in L_p^\ell}$, and 
$D^{p-1}(\A^H)$ is free with $\exp^{p-1}(\A^H)=\{e_J\}_{J \in L_{p-1}^{\ell-1}}$. Let $d:=|\A'|-|\A^H|$. 
Then $D^p(\A')$ has a free resolution
$$
0 
\rightarrow \oplus_{J \in L_{p-1}^{\ell-1}} S[-e_J-d-1]
\rightarrow \oplus_{J \in L_{p-1}^{\ell-1}} S[-e_J-d] \oplus (\oplus_{I 
\in L_p^\ell} S[-d_I] ) 
\rightarrow D^p(\A) \rightarrow 0.
$$
Here $$
L^{\ell}_p:=\{\{i_1,\ldots,i_p\} 
\mid 
1 \le i_1 <\cdots<i_p \le \ell\},
$$
and $$
d_I:=d_{i_1}+\cdots+d_{i_p}
$$
when $I=\{i_1,\ldots,i_p\}$.
\label{FR}
\end{cor}

The organization of this article is as follows. In \S2, we 
gather several results used for the proof of main results. 
In \S3 we introduce the free surjection argument that was introduced in \cite{A8} into the form that works in general contexts. \S4 is devoted to prove main results. 
In \S5 we introduce the dual version of the $B$-sequences. 
By using main results, we compute the 
Solomon-Terao polynomials of several arrangements, 
and give some answers to conjectures on $\Psi(\A;x,t)$ for some cases in 
\S6.
\medskip

\noindent
\textbf{Acknowledgements}.
The author is partially supported by JSPS KAKENHI Grant Numbers JP18KK0389 and JP21H00975.

\section{Preliminaries}

In this section let us introduce several definitions and results used for the proof of main results.

\begin{define}
    Let 
    $$
    L(\A):=\{
    \cap_{H \in \B} H\mid \B \subset \A\}
    $$
    be the \textbf{intersection lattice} of $\A$. 
The \textbf{M\"{o}bius function} $\mu:L(\A) 
\rightarrow \Z$ on $L(\A)$ is defined by $\mu(V)=1$ and by 
$$\mu(X)=
-\sum_{X \subsetneq Y \in L(\A)} \mu(Y).
$$
Then we can define the \textbf{characteristic polynomial} of $\A$ by 
$$
\chi(\A;t):=\sum_{X \in L(\A)}\mu(X)t^{\dim X}.
$$
Also for $X \in L(\A)$ the \textbf{localization} $\A_X$ of $\A$ at $X$ is defined by 
$$
\A_X:=\{H \in \A \mid X \subset H\},
$$
which is clearly a subarrangement of $\A$. Also the \textbf{restricion} 
$\A^X$ of $\A$ onto $X$ is defined by 
$$
\A^X:=\{X \cap H 
\mid H \in \A \setminus \A_X\},
$$
which is an arrangement in $X$.
\label{combin}
\end{define}

If we express 
$$
\chi(\A;t)=\sum_{i=0}^\ell (-1)^{\ell-i} b_i(\A) t^i,
$$
then $b_i(\A)$ is known to be the $i$-th Betti number of $M(\A):=
V \setminus \cup_{H \in \A} H$ when $\K=\CC$, see \cite{OS}. Thus 
$\chi(\A;t)$ is not only a combinatorial but also a topological invariant.

\begin{theorem}[The addition-deletion theorem, \cite{T1}]
Let $H \in \A$ and 
$\A':=\A \setminus \{H\}$.
Then two of the following three imply the third:
\begin{itemize}
\item[(1)]
$\A$ is free with $\exp(\A)=(d_1,d_2,\ldots,d_{\ell-1},d_\ell)$.
\item[(2)]
$\A'$ is free with $\exp(\A')=(d_1,d_2,\ldots,d_{\ell-1},d_\ell-1)$.
\item[(3)]
$\A^H$ is free with $\exp(\A^H)=(d_1,d_2,\ldots,d_{\ell-1})$.
\end{itemize}
Moreover, these three hold true if $\A$ and $\A'$ are 
free.
\label{addel}
\end{theorem}

\begin{define}[The Euler sequence, Proposition 4,45, \cite{OT}]
There is an exact sequence 
$$
0 \rightarrow D^p(\A') \stackrel{\cdot \alpha_H}{\rightarrow} D^p(\A) 
\stackrel{\rho^H}{\rightarrow} D^p(\A^H),
$$
Here $\rho^H(\theta)(\overline{f}):=\overline{\theta(f)}$, and 
$\overline{f} \in \overline{S}=S/\alpha_H S$ denotes the canonical image of $f \in S$ in $S/\alpha_H S$.
\label{Euler}
\end{define}

It is important to see when $\rho^H$ is surjective. One answer is as follows:

\begin{theorem}[Free surjection theorem (FST), Theorem 1.13, \cite{A9}]
If $\A'$ is free, then $\rho^H$ is surjective.
\label{FST}
\end{theorem}

The following is one of a few known results on $D^p(\A)$. 

\begin{theorem}[Proposition 3.4, \cite{ST}]
If $\A$ is free, then 
$D^p(\A)$ is also a free module. If $\theta_1,\ldots,\theta_\ell$ are 
basis for $D(\A)$, then 
$$
\{\theta_{i_1} \wedge \cdots \wedge \theta_{i_p}\}_{1 \le i_1 <\cdots <i_p \le \ell}
$$
form a basis for $D^p(\A)$.
\label{wedgefree}
\end{theorem}

The same result as Theorem \ref{wedgefree} exists for $\Omega^p(\A)$ when $\A$ is free. By using Theorem \ref{wedgefree} we can compute the Solomon-Terao polynomial of free arrangements as follows:

\begin{theorem}[\cite{ST}]
Assume that $\A$ is free with $\exp(\A)=(d_1,\ldots,d_\ell)$. Then 
$$
\Psi(\A;x,t)=\prod_{i=1}^\ell (-tx^{d_i}+1+\cdots+x^{d_i-1})
$$
and 
$$
\Psi(\A;x,-1)=\prod_{i=1}^\ell (1+\cdots+x^{d_i}).
$$
\label{STfree}
\end{theorem}

Next recall the classical polynomial $B$-theory, for the proof 
see Proposition 4.41 in \cite{OT}.

\begin{theorem}[\cite{T1}]
Let $H \in \A$ and $\A':=\A \setminus \{H\}$. Then there is $B \in S_{|\A'|-|\A^H|}$ such that 
$$
\theta (\alpha_H) \in (\alpha_H,B)
$$
for all $\theta \in D(\A')$. The polynomial $B$ is called 
\textbf{Terao's polynomial $B$}, and we can choose $B$ as $$
B=\displaystyle \frac{Q(\A')}{Q}.
$$
Here $Q \in S$ is an arbitrary lift up of $Q(\A^H) \in \overline{S}$ to $S$. 
\label{B}
\end{theorem}

Thus $B$ measures how far $\theta \in D(\A')$ is from $D(\A)$, and 
this led us to several important results on hyperplane arrangements like Theorem \ref{addel}, the division theorem in \cite{A2} and so on. A similar version for $\Omega^1(\A)$ is found in \cite{AD}. 
Since $D^0(\A)=S$, 
the exact sequence (\ref{seB}) when $p=1$ 
is nothing but Terao's polynomial $B$-theory. So the $B$-sequence (\ref{seB}) generalizes Terao's polynomial $B$ in Theorem \ref{B}.

%

%

Also let us recall a good structure for the logarithmic derivation module next to freeness. 

\begin{define}[\cite{A5}]
We say that $\A$ is \textbf{SPOG (strongly plus-one generated)} with 
$\POexp(\A)=(d_1,\ldots,d_\ell)$ and \textbf{level} $d$ if 
there is a minimal free resolution
$$
0 \rightarrow 
S[-d-1] \rightarrow \oplus_{i=1}^\ell S[-d_i] \oplus S[-d]
\rightarrow D(\A) \rightarrow 0.
$$
\label{SPOG}
\end{define}

The following shows the importance of SPOGness.

\begin{theorem}[Theorem 1.4, \cite{A5}]
Assume that $\A$ is free and $\A' :=\A \setminus \{H\}$ is not free for some $H \in \A$. Then $\A'$ is SPOG with $\POexp(\A')=\exp(\A)$ and 
level $|\A'|-|\A^H|$.
\label{SPOGthm}
\end{theorem}

\section{Free surjection theorem, revisited}

Free surjection theorem was first exhibited in \cite{A9}, and its origin is 
Theorems 1.8 and 1.13 in \cite{A8}. In this article we use the similar results to many different cases, so here we summarize it into the most general form. First let us introduce the local functor.

\begin{define}
Let $\A$ be an arrangement in $V=\K^\ell$ and 
let $F$ be a functor from arrangements in $V$ to modules over 
$S=\mbox{Sym}^*(V^*)$. We say that the functor $F$ is \textbf{local} if for any 
point $p \in \P^{\ell-1}=\mbox{Proj}(S)$, it holds that 
$$
F(\A_p)_p \simeq F(\A)_p,
$$
where $\A_p:=\{H \in \A \mid p \in H\}$. 
\label{localfunctor}
\end{define}

Now let us show the most general form of the surjection theorem.

\begin{theorem}[Projective dimensional surjection theorem]
Let $\A, \B$ be arrangements in $V=\K^\ell$, $H \in \A$ and 
let $\C$ an arrangement in $H \simeq \K^{\ell-1}$. 
Let $F$ be a local functor from arrangements in $V$ to modules over 
$S=\mbox{Sym}^*(V^*)$, with an exact sequence
$$
0 \rightarrow F(\B) \stackrel{f}{\rightarrow}
F(\A) \stackrel{g}{\rightarrow}
F(\C). 
$$
If $g$ is surjective at any localization $p \in H$ with 
$\codim p =k+2$, 
$\pd_S F(\B)=k< \ell-2$, and 
$$
H_*^0(\widetilde{F(\mathcal{D})})=F(\mathcal{D})
$$
for all arrangements $\mathcal{D}$,
then $g$ is surjective. 
\label{FST2}
\end{theorem}

\noindent
\textbf{Proof}. 
First let us prove that $g$ is locally surjective at any $p \in H$.
By the assumption it suffices to check for $p$ with $\codim p =s>k+2$. Let us prove by induction on $s
\ge k+2$. When $s=k+2$, the assumption completes the proof. So let us assume that $g$ is 
locally surjective along $H$ in codimension $s-1$. Thus we have 
$$
0 \rightarrow \widetilde{F(\B_p)} \stackrel{f}{\rightarrow}
\widetilde{F(\A_p)} \stackrel{g}{\rightarrow}
\widetilde{F(\C_p)} \rightarrow 0.
$$
Let $X=V(p)$. Then $\B_p$ decomposes into a direct product of 
an arrangement $\B^0$ in $X^{\perp} \simeq \K^s$ and the empty arrangement in $X$. 
So we may regard $\widetilde{F(B_p)}$ as a sheaf on $\P^{s-1}$. Now recalling that 
\begin{equation}
H^i_*(\widetilde{F(B)})=0\ (i=1,\ldots,\ell-k-2)
\label{eqpd}
\end{equation}
if $\B$ is an arrangement in $\K^\ell$ with a minimal free resolution
$$
0 \rightarrow F_{k} \rightarrow \cdots \rightarrow F_0 \rightarrow F(\B) \rightarrow 0,
$$
it holds that $H^1_*(\widetilde{F(\B^0)})=0$ since 
$$
s-k-2>k+3-k+2=1
$$
combined with (\ref{eqpd}), and an easy fact that 
$\pd_S F(\B_p) \le 
\pd_S F(\B)$ from the locality of $F$ and the exactness of the localization (see Theorem 2.2, \cite{KS} 
for example). Thus taking the global section, we know that 
$g$ is surjective along $H$ in codimension $s$. As a conclusion, $g$ is locally surjetive along $H$. Thus 
we have an exact sequence 
$$
0 \rightarrow \widetilde{F(\B)} \stackrel{f}{\rightarrow}
\widetilde{F(\A)} \stackrel{g}{\rightarrow}
\widetilde{F(\C)} \rightarrow 0.
$$
Now apply the same argument combined with (\ref{eqpd}) to complete the proof. \owari
\medskip

The following is an easy corollary of Theorem \ref{FST2} which is also 
often used.

\begin{theorem}[Free surjection theorem]
Let $\A, \B$ be arrangements in $V=\K^\ell$, $H \in \A$ and 
let $\C$ an arrangement in $H \simeq \K^{\ell-1}$. 
Let $F$ be a local functor from arrangements in $V$ to modules over 
$S=\mbox{Sym}^*(V^*)$, with an exact sequence
$$
0 \rightarrow F(\B) \stackrel{f}{\rightarrow}
F(\A) \stackrel{g}{\rightarrow}
F(\C). 
$$
Assume that $g$ is surjective at any localization $p \in H$, 
$F(\B)$ is free, and 
$$
H_*^0(\widetilde{F(\mathcal{D})})=F(\mathcal{D})
$$
for all arrangements $\mathcal{D}$,
then $g$ is surjective. Here 
$$
H_*^i(\widetilde{F(\mathcal{D})}):=\bigoplus_{j \in \Z} H^i(\P^{\ell-1},
\widetilde{F(\mathcal{D})}(j)).
$$
\label{FST}
\end{theorem}


The above properties hold for the logarithmic modules.

\begin{prop}[Proposition 6.6, \cite{ST}]
The functors 
\begin{eqnarray*}
\A \mapsto D^p(\A),\\
\A \mapsto \Omega^p(\A)
\end{eqnarray*}
are both local. 
\label{locality}
\end{prop}

\begin{prop}[Proposition 3.5, \cite{AD}]
$H^0_*(\widetilde{D^p(\A)})=D^p(\A)$ and 
$H^0_*(\widetilde{\Omega^p(\A)})=\Omega^p(\A)$ for 
all $0 \le p \le \ell$. 
\label{gs}
\end{prop}

\section{Proof of main results}

To prove Theorem \ref{higherB}, let us recall two results from \cite{AD} whose origin is in \cite{Z2}.

\begin{theorem}[Proposition 2.4, \cite{AD}]
Let $H \in \A$ and $\A':=\A \setminus \{H\}$. Then for each $0 \le p \le \ell$, 
there is an exact sequence 
$$
0 \rightarrow \Omega^{p}(\A) \stackrel{\cdot \alpha_H}{\rightarrow} 
\Omega^{p}(\A') \stackrel{i_H^p}{\rightarrow }\Omega^{p}(\A^H).
$$
Here $i^p(\omega)$ is defined as follows. For $\alpha_H=x_1$ and $\omega=\omega_0\wedge dx_1+\omega_1 
\in \Omega^p(\A')$, where $V'=\langle x_2,\ldots,x_\ell\rangle$ and 
$\omega_0,\omega_1 \in (\Omega^{p-1+i}_{V'})_{(0)}$, 
$$
i_H^p(\omega):=\overline{w_1}.
$$
\label{diffexact}
\end{theorem}

\begin{prop}[Lemma 2.3, \cite{AD}]
There is an identification 
$$Q(\A) \Omega^{p}(\A) \simeq D^{\ell-p}(\A)$$
by identifying $\wedge_{i \in I} dx_i$ with 
$\wedge_{i \not \in i} \partial_{x_i}$ for any $I \subset \{1,\ldots,\ell\}$ with 
$|I|=p$.
\label{ID}
\end{prop}

\noindent
\textbf{Proof of Theorem \ref{higherB}}.
By Theorem \ref{diffexact}, there is an exact sequence 
\begin{equation}
0 \rightarrow \Omega^{\ell-p}(\A) \stackrel{\cdot \alpha_H}{\rightarrow } 
\Omega^{\ell-p}(\A') \stackrel{i_H^{\ell-p}}{\rightarrow }\Omega^{\ell-p}(\A^H).
\label{eq16}
\end{equation}
Also recall the identification $Q(\A) \Omega^{\ell-p}(\A) \simeq D^p(\A)$ from Proposition \ref{ID}. 
Note the decomposition 
$Q(\A')=Q(\A^H)Q'$, and the fact that $\overline{Q'}=\overline{B}$ by 
definition of the polynomial $B$. 
Thus taking the product of the exact sequence (\ref{eq16}) with $Q(\A')$, we have 
$$
0 \rightarrow D^p(\A) \subset D^p(\A') 
\stackrel{\overline{\partial^p}}{\rightarrow} D^{p-1}(\A^H)\overline{B}.
$$
By the identification of $\wedge_{i \in I} dx_i $ with $\wedge_{i \not \in I} \partial_{x_i}$, and 
$\res_H^p(\wedge_{i \in I} dx_i)=\wedge_{i \in I} dx_i$ if $\alpha_H=x_1$ and $1 \not \in I$, and 
$0$ otherwise, 
we know that
$$
\overline{\partial}^p(\theta)(\overline{f_2},\ldots,\overline{f_p})=
\overline{\theta(\alpha_H,f_2,\ldots,f_p)},
$$
which completes the proof.\owari
\medskip

Thus the following is clear.

\begin{cor}
For $\theta \in D^p(\A')$, it holds that 
$$
\overline{\partial}^p(\theta) \in (\alpha_H,B)D_H^{p-1}(\A).
$$
Here 
$$
D_H^p(\A):=\{\theta \in D^p(\A) \mid \theta(\alpha_H,f_2,\ldots,f_p)=0\ (\forall f_2,\ldots,f_p \in S)\}.
$$
\label{higherB2}
\end{cor}

%

Now let us prove Theorem \ref{SPOGhigher}.
\medskip

\noindent
\textbf{Proof of Theorem \ref{SPOGhigher}}. 
Use the $B$-exact sequence. The statement is true if $\ell \le 2$. So assume that $\ell \ge 3$. 
Then Theorem \ref{FST2} shows that 
$\overline{\partial^p}$ is surjective, which completes the proof immediately by taking the Ext-exact sequence.\owari
\medskip

\noindent
\textbf{Proof of Corollary \ref{FR}}. 
By Theorems \ref{FST2} and \ref{FST}, we have the exact sequence 
\begin{equation}
0 \rightarrow D^p(\A) \rightarrow D^p(\A') \stackrel{\overline{\partial^p}}{\rightarrow} D^{p-1}(\A^H)\overline{B} \rightarrow 0.
\label{eq44}
\end{equation}
Let $\theta_1,\ldots,\theta_s$ be a basis for $D^p(\A)$ and 
$\varphi_1,\ldots,\varphi_t$ be the preimages of the basis for $D^{p-1}(\A^H)\overline{B}$ by $\overline{\partial^p}$. Here $s,t$ are appropriate combinatorial numbers. Let $\deg \theta_i=d_i$ and $\deg \varphi_i=e_i$. By (\ref{eq44}), it is clear that $\theta_1,\ldots,\theta_s,\varphi_1,\ldots,\varphi_t$ form a generator for $D^p(\A')$. So 
$$
0 \rightarrow K \rightarrow \oplus_i S[-d_i] \oplus \oplus_i S[-d-e_i] \rightarrow D^p(\A') \rightarrow 0
$$
holds. Let us check the relation in $K$. Since $\varphi_i \in D^p(\A') \setminus D^p(\A)$, there are relations 
$$
\alpha_H \varphi_i=\sum_{j=1}^s f_{ij} \theta_j
$$
at degree $d+e_i+1$ for $i=1,\ldots,t$. Let us show they generate $K$. Let 
$$
\sum_{i=1}^s g_i \theta_i+\sum_{i=1}^t h_i \varphi_i=0.
$$
Note that $\overline{\partial}^p(\theta_i)=0$. So if we send the above by $\overline{\partial^p}$, then all the coefficients $h_i$ has to be zero, i.e., 
$h_i$ can be replaced by $h_i \alpha_H$ and we have 
$$
\sum_{i=1}^s g_i \theta_i+\sum_{i=1}^t h_i (\alpha_H\varphi_i)=0.
$$
Since 
$$
\sum_{i=1}^t h_i (\alpha_H\varphi_i)=\sum_{i=1}^t \sum_{j=1}^s  h_if_{ij} \theta_j,
$$
the relation is transformed into the relation among $\theta_1,\ldots,\theta_s$. Namely, 
$$
(\sum_{i=1}^s g_i \theta_i+\sum_{i=1}^t h_i \alpha_H\varphi_i)
-(\sum_{i=1}^t h_i (\alpha_H \varphi_i-\sum_{j=1}^s f_{ij} \theta_j))
$$
has to be zero. So $K$ is generated by them. 
Since $\varphi_i$ appears exactly in the image of one basis element, the relations are 
independent. so we complete the proof. \owari
\medskip

\noindent
\textbf{Proof of Theorem \ref{STminus1dual}}. 
By Theorem \ref{FST2}, we have 
$$
0 \rightarrow D^p(\A') \stackrel{\cdot \alpha_H}{\rightarrow}
 D^p(\A) \stackrel{\rho^H}{\rightarrow}
D^p(\A^H) \rightarrow 0
$$
for all $p$. Thus 
$$
\mbox{Hilb}(D^p(\A);x)(t(x-1)-1)^p
=x\mbox{Hilb}(D^p(\A');x)(t(x-1)-1)^p+
\mbox{Hilb}(D^p(\A^H);x)(t(x-1)-1)^{p}.
$$
Now just take the sum by $p$ to conclude the proof.\owari
\medskip

\noindent
\textbf{Proof of Theorem \ref{STminus1}}. 
By Theorem \ref{FST2}, we have 
$$
0 \rightarrow D^p(\A) 
\rightarrow D^p(\A') \stackrel{\overline{\partial}^p}{\rightarrow}
D^{p-1}(\A^H)\overline{B} \rightarrow 0
$$
for all $p$. Thus 
\begin{eqnarray*}
\mbox{Hilb}(D^p(\A');x)(t(x-1)-1)^p
&=&\mbox{Hilb}(D^p(\A);x)(t(x-1)-1)^p\\
&+&
x^d (t(x-1)-1)\mbox{Hilb}(D^{p-1}(\A^H);x)(t(x-1)-1)^{p-1}.
\end{eqnarray*}
Now just take the sum by $p$ to conclude the proof.\owari
\medskip

\noindent
\textbf{Proof of Corollary \ref{STexplicit}}.
Clear by Theorems \ref{STminus1dual}, \ref{STminus1} and \ref{STfree}.
\medskip

Now recall 
the exact sequence (\ref{seB}) when $p=1$:
\begin{equation}
0 \rightarrow D(\A) \rightarrow 
D(\A') \stackrel{\overline{\partial}}{\rightarrow} \overline{S}\overline{B}.
\label{eq1}
\end{equation}

By using this we can prove, and understand several results naturally. 
For example, we can understand Theorem \ref{SPOGthm} naturally as follows:
\medskip


\noindent
\textbf{Another proof of Theorem \ref{SPOGthm}}. 
Assume that $\A'$ is not free. Since $\A$ is free, Theorem \ref{FST} shows that $\overline{\partial}^1$ is surjective. Thus there is $\theta \in D(\A')$ such that 
$\theta(\alpha_H)=B$ modulo $\alpha_H$. Hence 
$$
D(\A')=D(\A)+S\theta,
$$
which completes the proof. \owari
\medskip


\begin{rem}
Theorem \ref{SPOGthm} was proved by using Ziegler exact sequence, or the Euler exact sequence (\ref{euler}) in 
\cite{A5} and \cite{AD}, but they seem a bit strange from the viewpoint of the Euler exact sequence. Namely, if $\A$ is free and 
$\A'$ is not free, then 
the algebraic structure of $D(\A')$ is determined depending only on $|\A^H|$, independent of 
$D(\A^H)$. However, it is natural from the exact sequence (\ref{eq1}) because the third term is always $\overline{SB}$, which is free and determined by $\deg B=|\A'|-|\A^H|$!
\end{rem}

\section{Dual $B$-sequence}

Solomon-Terao formula is proved not only by using the logarithmic derivation modules but also the logarithmic differential forms as follows:

\begin{theorem}[Proposition 4.130, \cite{OT}]
Let
$$
\Phi(\A;x,t):=\sum_{p=0}^\ell \mbox{Hilb}(\Omega^p(\A);x)(t(1-x)-1)^p.
$$
Then 
$$
\chi(\A;t)=\Phi(\A;x,1).
$$
\label{STdual}
\end{theorem}

By Proposition \ref{ID}, 
we have the following.

\begin{prop}
$$
\Psi(\A;x,t)=x^{|\A|}(t(x-1)-1)^\ell\Phi(\A;x,\displaystyle \frac{-t}{t(x-1)-1}).
$$    
So
$$
\Psi(\A;x,-1)=(-1)^\ell x^{|\A|+\ell}\Phi(\A;x,\displaystyle \frac{-1}{x}).
$$    

\label{eq}
\end{prop}

\noindent
\textbf{Proof}. Clear by definition of these two polynomials.
\owari
\medskip

By using the completely same argument as for the logarithmic derivation modules, we can show the following $B$-sequence for the logarithmic differential forms.

\begin{theorem}[Dual $B$-sequence]
Let $H \in \A$ and  $\A':=\A \setminus \{H\}$.
Then there is an exact sequence
$$
0 \rightarrow \Omega^p(\A') \rightarrow 
\Omega^p(\A) \stackrel{\res^p}{\rightarrow}
\Omega^{p-1}(\A^H)/\overline{B}.
$$
Here $B$ is  the polynomial $B$ appeared in Theorem \ref{B}, and for $x_1=\alpha_H$ and the decomposition 
$$
\Omega^p(\A) \ni \omega=\omega_0\wedge \displaystyle \frac{dx_1}{x_1}
+\omega_1
$$
for two logarithmic differential forms $\omega_0,\omega_1 \in (\Omega^{p-1+i}_{V'})_{(0)}
$ with $V'=\langle x_2,\ldots,x_\ell \rangle$, $\res^p$ is defined by 
$$
\res^p(\omega):=\overline{\omega}_0.
$$
\label{dualB}
\end{theorem}

\noindent
\textbf{Proof}.
The same as the proof of Theorem \ref{seB}.\owari
\medskip

Thus we have the same structure theorem for higher logarithmic differential forms as 
follows:

\begin{theorem}
Let $H \in \A$ and  $\A':=\A \setminus \{H\}$.
Assume that $\Omega^p(\A')$ is free. If 
$\Omega^p(\A)$ is not free, then 
$$
\Omega^p(\A)=\Omega^p(\A')+\sum_{i=1}^s S \omega_i,
$$
where $\res^p(\omega_1),\ldots,\res^p(\omega_s)$ form a minimal set of generators for $\Omega^{p-1}(\A^H)/\overline{B}$. Moreover, 
$$
\pd_S \Omega^p(\A) =\pd_{\overline{S}} \Omega^{p-1}(\A^H)+1.
$$
In particular, if $\Omega^p(\A')$ and $\Omega^{p-1}(\A^H)$ are both free, then $\pd_S \Omega^p(\A) \le 1$. 
\label{SPOGhigherdual}
\end{theorem}



\begin{cor}
Assume that $\Omega^p(\A')$ is free with $\exp^p(\A)=\{d_I\}_{I \in 
L_p^\ell}$, and 
$\Omega^{p-1}(\A^H)$ is free with $\exp^{p-1}(\A^H)=\{e_J
\}_{J \in L_{p-1}^{\ell-1}}$. Let $d:=|\A'|-|\A^H|$. 
Then $\Omega^p(\A)$ has a free resolution
$$
0 
\rightarrow \oplus_{J \in L^{\ell-1}_{p-1} }S[e_J+d]
\rightarrow \oplus_{J \in L_{p-1}^{\ell-1}} S[e_J+d+1] \oplus (\oplus_{
I \in L^\ell_p
} S[d_I] ) 
\rightarrow \Omega^p(\A) \rightarrow 0.
$$
\label{FRdual}
\end{cor}

\noindent
\textbf{Proof}. The same proof as that of Corollary \ref{FR}. Note that 
the residue map in the dual $B$-sequence increases the degree by one. \owari
\medskip


\begin{rem}
We can show Theorems \ref{STminus1dual} and \ref{STminus1} by using 
the dual $B$-sequences.
\end{rem}

\section{Examples and conjectures}

As mentioned in \S1, to compute $\Psi(\A;x,t)$ is really hard. The only known case is when $\A$ is free with $\exp(\A)=(d_1,\ldots,d_\ell)$ as seen in Theorem \ref{STfree}.
Now by Theorems \ref{STminus1dual} and \ref{STminus1}, we have tools to compute the Solomon-Terao 
polynomial for non-free arrangements. Let us check some of them.

\begin{example}
Let 
$$
\A:=xyzw(x+y)(x+z)(x+w)(x+y+z)(x+y+w)(x+z+w).
$$
$\A$ is free with $\exp(\A)=(1,3,3,3)$ and let $H:\{y=0\} \in \A$. Then $\A^H$ is free with $\exp(\A^H)=(1,2,3)$. So by the 
deletion theorem and Theorem \ref{SPOGthm}, $\A':= \A \setminus \{H\}$ is not free but 
SPOG with $\POexp(\A')=(1,3,3,3)$ and level $3$. 
Now we can apply Theorem \ref{STminus1} to compute 
\begin{eqnarray*}
\Psi(\A';x,t)&=&
\Psi(\A;x,t)-x^{3} (t(x-1)-1)\Psi(\A^H;x,t)\\
&=&
(tx-1)(tx^3-1-x-x^2)^3\\
&\ &-
x^3(t(x-1)-1) (tx-1)(tx^2-1-x)(tx^3-1-x-x^2)\\
&=&
(1+3x+6x^2+6x^3+4x^4+x^5)\\
&\ &-(x+3x^2+10x^3+13x^4+12x^5+4x^6)t\\
&\ &+(4x^4+8x^5+12x^6+6x^7)t^2\\
&\ &-(x^6+4x^7+4x^8)t^3\\
&\ &+x^9t^4.
\end{eqnarray*}
Hence 
\begin{eqnarray*}
\Psi(\A';x,-1)&=&
1+4x+9x^2+16x^3+21x^4+21x^5+17x^6+10x^7+4x^8+x^9.
\end{eqnarray*}
\label{X3}
\end{example}

\begin{example}
By using the same $\A$ as in Example \ref{X3}, 
let us consider the Solomon-Terao polynomial 
that is a bit more far from free arrangements. i.e., 
let $\A$ be an arrangement in $\R^4$ defined by 
$$
Q(\A)=xyzw(x+y)(x+z)(x+w)(x+y+z)(x+y+w)(x+z+w).
$$
Contrary to Example \ref{X3}, let $H:x=0$ and $\A':=\A \setminus \{H\}$. Then Theorem \ref{STminus1} shows that 
$$
\Psi(\A';x,t)=\Psi(\A;x,t)-x^3(t(x-1)-1)\Psi(\A^H;x,t).
$$
Since $\A$ is free with $\exp(\A)=(1,3,3,3)$, we know that 
$$
\Psi(\A;x,t)=(-tx+1)(-tx^3+1+x+x^2)^3.
$$
So we need to compute $\Phi(\A^H;x,t)$. However, $\A^H$ is not free but SPOG 
with $\POexp(\A^H)=(1,3,3)$ and level $3$. In fact, this $\A$ and $\A^H$ is the smallest 
counter example to Orlik's conjecture asserting that every restriction of a free arrangement is free. 
Let 
$$
Q(\A^H)=yzw(y+z)(y+w)(z+w)
$$
and let $X:z+w=0$ in $H$. Then it is easy to show that 
$\A^H \setminus \{X\}$ is free with exponents $(1,2,2)$. Since $|(\A^H)^X|=4$, we can apply 
Theorem \ref{STminus1dual} to compute 
\begin{eqnarray*}
\Psi(\A^H;x,t)&=&x\Psi(\A^H \setminus \{X\};x,t)+\Psi(\A^X;x,t)\\
&=&x(-tx+1)(-tx^2+1+x)^2+(-tx+1)(-tx^3+1+x+x^2).
\end{eqnarray*}
Combining these, we obtain
\begin{eqnarray*}
\Psi(\A';x,t)&=&\Psi(\A;x,t)+x^3(t(x-1)-1)\Psi(\A^H;x,t)\\
&=&(-tx+1)((-tx^3+1+x+x^2)^3\\
&\ &+x^3(t(x-1)-1)((-tx^2+1+x)^2+(-tx^3+1+x+x^2))\\
&=&(-tx+1)(-t^3x^8+t^2(2x^7+6x^6)-t(x^6+10x^5+7x^4+4x^3)\\
&\ &+x^6+x^5+3x^4+5x^3+6x^2+3x+1).
\end{eqnarray*}
Thus 
\begin{eqnarray*}
\Psi(\A';x,-1)&=&
(1+x)(1+3x+6x^2+9x^3+10x^4+11x^5+8x^6+2x^7+x^8).
\end{eqnarray*}
\label{X3-2}
\end{example}

The argument used in Example \ref{X3-2} is generalized as follows:

\begin{theorem}
Let $H \in \A,\ \A':=\A \setminus \{H\}$. Let $d:=|\A'|-
|\A^H|$.  
\begin{itemize}
    \item[(1)]
    Assume that $\A$ is free, and there is $X \in \A^H$ such that 
    $\B:=\A^H \setminus \{X\}$ is free. Then 
    $$
    \Psi(\A';x,t)=\Psi(\A;x,t)+x^d(t(x-1)-1)(x\Psi(\B;x,t)+\Psi(\A^X;x,t)).
    $$
    \item[(2)]
    Assume that $\A$ is free, and there is $X \not \in \A^H$ such that 
    $\B:=\A^H \cup \{X\}$ is free. Then 
    $$
    \Psi(\A';x,t)=\Psi(\A;x,t)+x^d(t(x-1)-1)(\Psi(\B;x,t)+x^e(t(x-1)-1)\Psi(\A^X;x,t)).
    $$
    Here $e:=|\A^H|-1-|\A^X|$.
    \item[(3)]
    Assume that $\A'$ is free, and there is $X \in \A^H$ such that 
    $\B:=\A^H \setminus \{X\}$ is free. Then 
    $$
    \Psi(\A;x,t)=x\Psi(\A';x,t)+x\Psi(\B;x,t)+\Psi(\A^X;x,t).
    $$
    \item[(4)]
    Assume that $\A'$ is free, and there is $X \not \in \A^H$ such that 
    $\B:=\A^H \cup \{X\}$ is free. Then 
    $$
    \Psi(\A;x,t)=x\Psi(\A';x,t)+\Psi(\B;x,t)+x^e(t(x-1)-1)\Psi(\A^X;x,t)).
    $$
    Here $e:=|\A^H|-1-|\A^X|$.
\end{itemize}
\end{theorem}

\noindent
\textbf{Proof}.
Apply Theorems \ref{STminus1dual} and \ref{STminus1} continuously. \owari
\medskip

\begin{rem}
Note that all the freeness of arrangements appearing in Example \ref{X3} are divisionally free, see 
\cite{A2}. Thus their freeness depends only on $L(\A)$, i.e., combinatorially determined. So 
$\Psi(\A';x,t)$ in Example \ref{X3} is combinatorial. 

In general, 
if an arrangement $\A$ can be constructed by adding/deleting hyperplanes in such a way that all the appearing 
arrangements in that process are combinatorially free (like divisionally free in \cite{A2} or 
additively free in \cite{A6}), then its Solomon-Terao polynomual is combinatorially determined by using Theorems \ref{STminus1dual} and \ref{STminus1}.
\end{rem}

\begin{example}
    Let $\A:=xyz(x+y+z)$,\ 
    $H=\{x+y+z=0\}$ and let 
    $\A'=\A \setminus\{H\}$. Then $\A'$ is free with $\exp(\A')=(1,1,1)$ and 
    $\exp(\A^H)=(1,2)$, both free. 
    In Example 5.4 in \cite{AMMN}, $\Psi(\A;x,-1)$ was exhibited 
    to be the case when it is hard to compute even though $\A'$ and $\A^H$ are both free. Now we can use Theorem \ref{STminus1dual} to compute 
    \begin{eqnarray*}
\Psi(\A;x,-1)&=&x\Psi(\A';x,-1)+\Psi(\A^H;x,-1)\\
&=&x(x+1)^3+(x+1)(1+x+x^2)\\
&=&1+3x+5x^2+4x^3+x^4.
    \end{eqnarray*}
Moreover, by Theorem \ref{STminus1dual} again, 
\begin{eqnarray*}
\Psi(\A;x,t)&=&x\Psi(\A';x,t)+\Psi(\A^H;x,t)\\
&=&x(-tx+1)^3+(-tx+1)(-tx^2+1+x)\\
&=&2x+1-t(5x^2+x)+4x^3t^2-x^4t^3.
    \end{eqnarray*}
\label{3generic}
\end{example}

Example \ref{3generic} can be generalized in the following manner:

\begin{theorem}
    Let $\A$ be a generic arrangement in $\K^\ell$ consisting of $(\ell+1)$-hyperplanes. Then 
    \begin{eqnarray*}
    \Psi(\A;x,t)&=&
    (1+(\ell-1)x)-t(x+(\begin{pmatrix}\ell+1\\2\end{pmatrix}-1)x^2)\\
    &+&\sum_{i=2}^\ell (-1)^i \ \begin{pmatrix}\ell+1\\i+1\end{pmatrix}x^{i+1} t^i.
    \end{eqnarray*}
    \label{generic}
\end{theorem}

\noindent
\textbf{Proof}.
Apply the induction on $\ell$ starting from Example 
\ref{3generic} with Theorem \ref{STexplicit} (1). \owari
\medskip

Related to the Solomon-Terao polynomials, there are many open problems posed in \cite{AMMN}. Let us pose two of them.

\begin{conj}[Conjecture 5.12, \cite{AMMN}]
$\deg \Psi(\A;x,-1)=|\A|$ and it is monic.
\label{degree}
\end{conj}

\begin{conj}[Conjecture 3.5, \cite{AMMN}]
$\A$ is free if and only if $\Psi(\A;x,-1)$ is palindromic, or equivalently, 
$$
\Psi(\A;x,-1)=\prod_{i=1}^\ell \displaystyle \frac{1-x^{d_i+1}}{1-x}
$$
for some $d_1,\ldots,d_\ell$.
\label{palindromic}
\end{conj}

By the above theorems, we can give a partial answer to these conjectures.

\begin{cor}
Let $\ell \ge 3$, and assume that all the arrangements appearing here are essential. 
\begin{itemize}

    \item [(1)]
If $\A$ and $\A^H$ are free, then Conjectures \ref{degree} and \ref{palindromic} are true for $\A':=\A \setminus \{H\}$ for any $H \in \A$.
\item [(2)]
If $\A$ and $\A^H$ are free, then Conjectures \ref{degree} and \ref{palindromic} are true for $\A:=\A' \cup \{H\}$ for any $H \not \in \A$.
    \end{itemize}
    \label{ell3}
\end{cor}

\noindent
\textbf{Proof}.
(1)\,\,
Use the formula in Theorem \ref{STminus1}:
\begin{equation}
    \label{formula}
\Psi(\A';x,-1)=
\prod_{i=1}^\ell 
(1+x+\cdots+x^{d_i})-x^{d+1}\prod_{i=1}^{\ell -1}
(1+x+\cdots+x^{e_i}).
\end{equation}
Here we use the same notations of exponents as in Theorem \ref{STminus1}. 
The highest degree term in the first product in (\ref{formula}) is $x^{|\A|}$ and 
that of the second product in (\ref{formula}) is $-x^{|\A|}$ since $$
d+1+\sum_{i=1}^{\ell-1}e_i=
|\A'|-|\A^H|+1+|\A^H|=|\A|.
$$
Note that $1=d_1=e_1$ and $2 \le d_2,\ 2 \le e_2$ since $\A$ and $\A^H$ are essential. 
Thus the second highest degree term in the first product in (\ref{formula}) is $$
(1+\ell-1)x^{|\A|-1}=\ell x^{|\A|-1}.
$$ 
On the other hand, 
the second highest degree term in the second product in (\ref{formula}) is $$
-(1+\ell-2)x^{|\A|-1}=-(\ell-1) x^{|\A|-1}.
$$ 
Thus $\Psi(\A';x,-1)$ is a monic polynomial of degree 
$|\A'|$, which confirms Conjecture \ref{degree} in this case. 

Next let us check Conjecture \ref{palindromic} for (1). 
It is clear that $d \ge 2$. 
By (\ref{formula}), the coefficients of $x$ and $x^{|\A|-2}$ in $\Psi(\A';x,-1)$ are
$\ell$. So compute the coefficient of $x^{2}$, that is 
$$
(\ell-1)+\ell(\ell-1)/2,
$$
and that of $x^{|\A|-3}$ is 
\begin{eqnarray*}
&\ &(\ell-1)+\ell(\ell-1)(\ell-2)/6+
(\ell-1)(\ell-2)/2+(\ell-n_2(\A^H))\\
&-&
(
(\ell-2)+(\ell-1)(\ell-2)(\ell-3)/6+
(\ell-2)(\ell-3)/2+(\ell-1-n_2(\A)).
)\\
&=&(\ell-1)+(\ell-1)(\ell-2)/2+n_2(\A^H)-n_2(\A).
\end{eqnarray*}
Here $n_2(\A)$ denotes the cardinality of $2$ appearing in $\exp(\A)$. So if these two are the same, then it has to hold that 
$$
n_2(\A^H)-n_2(\A)=\ell-1.
$$
Since the degree of the Euler derivation, $1$, is contained in $\exp(\A^H)$, this is impossible. 
So $\Psi(\A';x,-1)$ is not palindromic, confirming Conjecture \ref{palindromic} in this case.

(2)\,\,
Apply the same argument as in (1) combined with Theorem \ref{STminus1dual}. \owari


\begin{thebibliography}{ABCHT}




\bibitem{A2}
T. Abe,
Divisionally free arrangements of hyperplanes. 
\textit{Invent. Math.} \textbf{204} (2016), no. 1, 317--346.






\bibitem{A5}
T. Abe, 
Plus-one generated and next to free arrangements of hyperplanes. 
\textit{Int. Math. Res. Not} \textbf{2021} (2021), no. 12, 9233--9261.




\bibitem{A8}
T. Abe, 
Double points of free projective line arrangements, 
\textit{Int. Math. Res. Not.} 
\textbf{2022} (2022), no. 3, 1811--1824.

\bibitem{A6}
T. Abe, 
Addition-deletion theorem for free hyperplane arrangements and 
combinatorics.
\textit{J. Algebra} \textbf{610} (2022), 1--17.

\bibitem{A9}
T. Abe, 
Projective dimensions of hyperplane arrangements.
arXiv:2009:04101 (2020).












\bibitem{AD}
T. Abe and G. Denham, 
Deletion-restriction for logarithmic forms on multiarrangements. 
arXiv:2203.04816 (2022).





\bibitem{AHMMS}
T. Abe, T. Horiguchi,  M. Masuda, S. Murai and T. Sato, 
Hessenberg varieties and hyperplane arrangements. 
\textit{J. Reine Angew. Math.} 
\textbf{764} (2020), 241--286.

\bibitem{AMMN}
T. Abe, T. Maeno, S. Murai 
and Y. Numata, 
Solomon-Terao algebra for hyperplane arrangements.
\textit{J. Math. Soc. Japan} \textbf{71} (2019), no. 4, 1027--1047.
































 \bibitem{KS}
J. Kung and H. Schenck, 
Derivation modules of orthogonal duals of hyperplane arrangements.
\textit{J. Algebraic Combin}. 
\textbf{24} (2006), no. 3, 253--262.
 

 \bibitem{OS}
 P. Orlik and L. Solomon, 
Combinatorics and topology of complements of hyperplanes. 
\textit{Invent. Math}. \textbf{56} (1980), no. 2, 167--189.

\bibitem{OT} P. Orlik and H. Terao, \textit{Arrangements of hyperplanes}.
Grundlehren der Mathematischen Wissenschaften, 
\textbf{300}. Springer-Verlag, Berlin, 1992.




\bibitem{Sa}
K. Saito, 
Theory of logarithmic differential forms and logarithmic vector fields.
\textit{J. Fac. Sci. Univ. Tokyo} \textbf{27} (1980), 265--291.   





\bibitem{ST}
L. Solomon and H. Terao, 
A formula for the characteristic polynomial
of an arrangement. \textit{Adv. Math.} \textbf{64} (1987), 
305--325.

\bibitem{T1}
H. Terao, 
Arrangements of hyperplanes and their freeness I, II. 
\textit{J. Fac. Sci. Univ. Tokyo} \textbf{27} (1980), 293--320.   

\bibitem{T2}
H. Terao, 
Generalized exponents of a free arrangement of hyperplanes and
Shephard-Todd-Brieskorn formula. \textit{Invent. math}. 
\textbf{63}  (1981),
159--179.




\bibitem{Y1}
M. Yoshinaga,
Characterization of a free arrangement and
conjecture of
Edelman and Reiner. \textit{Invent. Math.} \textbf{157} (2004), no. 2,
449--454.






\bibitem{Z2}
G. M. Ziegler, Combinatorial construction of logarithmic differential forms.
\textit{Adv. Math}. \textbf{76} (1989), 116--154.%
	
\end{thebibliography}
\end{document}